\documentclass{article}
\usepackage{amsmath,amsthm,amssymb,amsfonts,latexsym,graphicx}
\usepackage{makeidx}
\usepackage{bbm}
\usepackage{mathabx}
\usepackage[margin=1in]{geometry}
\usepackage[shortlabels]{enumitem}
\usepackage{theoremref}
\usepackage{authblk}

\newcommand{\N}{{\mathbb{N}}}
\newcommand{\R}{{\mathbb{R}}}
\newcommand{\Z}{{\mathbb{Z}}}
\newcommand{\C}{{\mathbb{C}}}

\DeclareMathOperator\supp{supp}

\newtheorem{theorem}{Theorem}
\newtheorem{lemma}{Lemma}

\title{$L^p$ Estimates for Eigenfunctions of a Generalized Landau Magnetic Laplacian}
\author{Ben Goldschlager}
\affil{University of Maryland, College Park}
\date{}

\begin{document}

\maketitle

Email: bgoldsc1@umd.edu

ORCiD: 0009-0007-5365-1623

\newpage

\begin{abstract}
In this paper, we examine eigenfunctions of a generalized Landau Magnetic Laplacian that models the physics of an electron confined to a plane in a magnetic field orthogonal to the plane. This operator has an infinite dimensional null space and, at least in the model case, has infinite dimensional eigenspaces with eigenvalues which are essentially the same as the eigenvalues of the Hermite operator. We demonstrate that, under fairly general assumptions on the potential function of the magnetic field, the $L^\infty$ norm of these eigenfunctions is bounded by their $L^2$ norm independently of the associated eigenvalue. We furthermore demonstrate an improvement in the $L^6$ norm of these eigenfunctions. The method we use comes from semiclassical analysis and is inspired by the work of Koch, Tataru, and Zworski that applies locally. In our case, we use a new conjugation argument to demonstrate the result over all of $\R^2$.    
\end{abstract}


\section{Introduction}

In this paper, we examine eigenfunctions of a generalized Landau magnetic Laplacian, $H = D^* D$, where $D = \partial_z + \partial_z \phi, z \in \C$ with real-valued scalar potential function $\phi$. This operator arises as the Hamiltonian that describes an electron confined to a plane in a magnetic field orthogonal to that plane, see \cite{XD} and \cite{Fefferman} for further background. In a suitable gauge, the vector potential, $A$, of the magnetic field is

\begin{equation*}
    A = \begin{bmatrix} - \partial_2 \phi \\ \partial_1 \phi \\ 0 \end{bmatrix}
\end{equation*}

Our goal is as follows. We suppose

\begin{equation}
    |\partial^\alpha \phi| \leq C_\alpha \text{ for all } |\alpha| \geq 2
\end{equation}

and consider eigenfunctions $u$ of $H$ where

\begin{equation} \label{geigenfunction_eq}
    Hu = \lambda^2 u
\end{equation}

Then, we will demonstrate the following two results.

\begin{theorem} \thlabel{main}
    With $H$ and $\phi$ defined as above, we have

    \begin{equation}
        ||u||_{L^\infty} \lesssim ||u||_{L^2}
    \end{equation}

    where the implicit constant is independent of the eigenvalue $\lambda^2$.
\end{theorem}

and

\begin{theorem} \thlabel{main_L6}
    With $H$ and $\phi$ defined as above, we have

    \begin{equation}
        ||u||_{L^6} \lesssim \lambda^{-1/3}||u||_{L^2}
    \end{equation}
\end{theorem}

Note that the second result is strictly better than what we would get by interpolation with the $L^\infty$ estimate. 

In the specific case that $\phi(z) = |z|^2$, \thref{main} is contained implicitly in the work of Folland and Thangavelu, see \cite{Folland} and \cite{Thangavelu}. Folland, for example, shows that functions in the null space of this operator have the form $e^{-|z|^2} F(\bar{z})$ where $F$ is anti-holomorphic. Then, the null space has the reproducing kernel $K(z, w) = e^{2z \bar{w}}$, from which the $L^\infty$ estimate easily follows. He furthermore shows that the eigenspaces have eigenvalues which are the even non-negative integers (called Landau levels) and that functions in these eigenspaces can be obtained by applying the corresponding creation operator to eigenfunctions in the lower eigenspaces. Thus, the eigenspaces are infinite dimensional. Furthermore, Thangavelu's work demonstrates that these eigenfunctions can be written as twisted convolutions with a special Hermite function, in which case the $L^\infty$ result follows from a simple convolution inequality and growth estimates on the special Hermite functions. To the best of my knowledge, however, no previous work has been done on the $L^6$ inequality. 

In the case of a general potential function, it remains true that the null space contains functions of the form $e^{-\phi} F(\bar{z})$ where $F$ is anti-holomorphic, provided $\phi$ grows fast enough that these functions are $L^2$. For example, in the case that $\phi \gtrsim |z |^\epsilon, \epsilon > 0$, the null space will contain functions of the form $e^{- \phi} p(\bar{z})$ where $p$ is any polynomial. Thus, the null space of $H$ will still be infinite dimensional. We furthermore expect that, in the case of compact perturbations of the original operator with potential function $\phi(z) = |z|^2$, the higher eigenspaces may also be infinite dimensional or may split into spectral clusters.

Note that, by switching to real coordinates, we can rewrite $H$ as

\begin{equation}
    H = \left(\frac{D_1}{2} - \frac{\partial_2 \phi}{2}\right)^2 + \left(\frac{D_2}{2} + \frac{\partial_1 \phi}{2}\right)^2 - \frac{\Delta \phi}{4}
\end{equation}
 
which is the form of the operator with which we will concern ourselves hereafter. 

For example, in the specific case that $\phi(z) = |z|^2$, this becomes

\begin{equation}
    H = \left(\frac{D_1}{2} - x_2\right)^2 + \left(\frac{D_2}{2} + x_1 \right)^2 - 1
\end{equation}

Thus, we see that, when $x$ is confined to a compact ball, this operator behaves like the standard Laplacian, and we describe below how we deal with the case of $x$ large using a conjugation argument. 

In order to prove these results, we take advantage of developments in semiclassical analysis; see \cite{ktz} and \cite{Zworski}. In particular, our method of proof is inspired by work done by Koch, Tataru, and Zworski in \cite{ktz}, which generalizes earlier work of Koch and Tataru on the Hermite operator, see \cite{kochtataru}, and of Sogge on the Laplace-Beltrami operator on a compact Riemannian manifold, see \cite{Sogge}. Like the Laplace-Beltrami operator, the operator we consider is also of principal type. However, it has the same scaling as in the case of the Hermite operator. Koch, Tataru, and Zworski are able to prove estimates on eigenfunctions of principal type operators by factoring into an elliptic operator and a Schrodinger operator and using energy estimates and Strichartz estimates for the Schrodinger operator. In our case, we likewise factor our operator into an elliptic operator and a Schrodinger operator and use the same energy and Strichartz estimates that Koch, Tataru, and Zworski do. We then arrive at results that are analogous to the results of Sogge for the Laplace-Beltrami operator on a compact Riemannian manifold. 

What sets our approach apart is that we are working on all of $\R^2$ and therefore some care is needed to appropriately localize our eigenfunctions. We take advantage of a multiplication operator which we denote $T_q, q \in \R^2$, such that conjugating our operator with $T_q$ allows us to approximately translate our operator in a manner which is uniformly bounded in $q$. 

This paper is organized as follows. In section 2, we prove \thref{main} and \thref{main_L6} for the case when eigenvalues $\lambda^2$ are sufficiently large. In section 3, we prove \thref{main} and \thref{main_L6} for eigenvalues $\lambda^2$ which are less than some constant. We have an appendix at the end which contains some important results of semiclassical analysis which we use in the body of this paper. 

\section{Large Eigenvalue Case}

In this section, we turn to the large eigenvalue case, which is the more difficult and interesting case.

Now, if we let $u_\lambda (x) = u(\lambda x)$ and $\phi_\lambda (x) = \phi(\lambda x)$, then \eqref{geigenfunction_eq} gives us

\[\left[\left(\frac{1}{2} \frac{D_1}{\lambda} - \frac{\partial_2 \phi_\lambda}{2 \lambda}\right)^2 + \left(\frac{1}{2} \frac{D_2}{\lambda} + \frac{\partial_1 \phi_\lambda}{2 \lambda}\right)^2 - \frac{\Delta \phi_\lambda}{4 \lambda^2}\right] u_\lambda = \lambda^2 u_\lambda\]

Dividing by $\lambda^2$ and rearranging, we get

\[\left[\left(\frac{1}{2}\frac{D_1}{\lambda^2} - \frac{\partial_2 \phi_\lambda}{2 \lambda^2}\right)^2 + \left(\frac{1}{2} \frac{D_2}{\lambda^2} + \frac{\partial_1 \phi_\lambda}{2 \lambda^2}\right)^2 - \frac{\Delta \phi_\lambda}{4 \lambda^4} - 1 \right] u_\lambda = 0\]

If $h = \frac{1}{\lambda^2},$ let $u_h(x) = u(h^{-1/2}x) = u(\lambda x) = u_\lambda(x)$ and $\phi_h(x) = \phi(h^{-1/2}x) = \phi(\lambda x) = \phi_\lambda(x)$, then this becomes

\[\left[\left(\frac{h}{2} D_1 - \frac{h}{2} \partial_2 \phi_h\right)^2 + \left(\frac{h}{2} D_2 + \frac{h}{2} \partial_1 \phi_h\right)^2 - \frac{h^2}{4} \Delta \phi_h - 1\right] u_h = 0\]

Or, by rewriting this, we have

\[\left[\left(\frac{h}{2} D_1 - \frac{h^{1/2}}{2} (\partial_2 \phi)(h^{-1/2}x)\right)^2 + \left(\frac{h}{2} D_2 + \frac{h^{1/2}}{2} (\partial_1 \phi)(h^{-1/2}x)\right)^2 - \frac{h}{4} (\Delta \phi)(h^{-1/2}x) - 1\right] u_h =: Pu_h = 0\]

Note that $||u_h||_{L^\infty} = ||u||_{L^\infty}$ and, by a change of variables, $||u_h||_{L^2} = h^{1/2}||u||_{L^2}$ and $||u_h||_{L^6} = h^{1/6} ||u||_{L^6}$.

Thus, \thref{main} can be reformulated as follows.

\begin{theorem} Suppose $P u_h = 0$, then

\begin{equation}
    ||u_h||_{L^\infty} \lesssim h^{-1/2} ||u_h||_{L^2}
\end{equation}

\end{theorem}

Likewise, \thref{main_L6} can be reformulated as

\begin{theorem} Suppose $P u_h = 0$, then

\begin{equation}
    ||u_h||_{L^6} \lesssim h^{-1/6} ||u_h||_{L^2}
\end{equation}

\end{theorem}

Define 

\begin{align*}
    A & := \frac{h}{2} D_1 - \frac{h}{2} (\partial_2 \phi_h)(x) \\ 
    B & := \frac{h}{2} D_2 + \frac{h}{2} (\partial_1 \phi_h)(x) \\
    P & := A^2 + B^2 - \frac{h^2}{4} (\Delta \phi_h)(x) - 1 \\
    \tilde{A}_q &:= \frac{h}{2} D_1 - \frac{h}{2} (\partial_2 \phi_h)(x+q) + \frac{h}{2} (\partial_2 \phi_h)(q)\\
    \tilde{B}_q & := \frac{h}{2} D_2 + \frac{h}{2} (\partial_1 \phi_h)(x+q) - \frac{h}{2} (\partial_1 \phi_h)(q) \\
    \tilde{P}_q & :=  \tilde{A}_q^2 + \tilde{B}_q^2 - \frac{h^2}{4} (\Delta \phi_h)(x) - 1 \\
\end{align*}

Likewise, we use the corresponding lowercase letters to denote the symbols corresponding to the above operators, e.g. $A = a^\mathrm{w}(x, hD)$ and $\tilde{A}_q = \tilde{a}_q^\mathrm{w}(x, hD)$. Note, furthermore, that, in the specific case of a quadratic potential function, $\tilde{A}_q = A$ and likewise for $\tilde{B}_q$ and $\tilde{P}_q$.

We likewise define

\begin{align*}
    u_{h,q}(x) &:= u_h(x-q) \\
    \sigma(x, \xi) &:= \langle x_2, \xi_1 \rangle - \langle x_1, \xi_2 \rangle \\
    T_q &:= e^{i \sigma(x, h^{-1/2}\nabla \phi(h^{-1/2}q)) } \\
    S_{\rho, \delta}(m) &:= \{ a \in C^\infty | |\partial_\xi^\alpha \partial_x^\beta a| \lesssim h^{-\rho |\alpha| - \delta |\beta|} m \text{ for all multiindices } \alpha, \beta \} 
\end{align*}

Here, $\sigma$ is the standard symplectic product. $T_q$ allows us to approximately translate our operators in a manner which is uniform in $q$. The uniformity in $q$ follows from our assumptions on $\phi$. For example, we see that 

\begin{equation*}
    T_q^{-1} \tilde{A}_q T_q = \frac{h}{2} D_1 - \frac{h}{2} (\partial_2 \phi_h)(x + q)
\end{equation*}

Finally, $S_{\rho, \delta}(m)$ is an extension of the symbol classes in \cite{Zworski}, about which we say more in the appendix.

Let $\beta \in C_0^\infty(\R^2)$ such that $\beta \equiv 1$ for $|x| \leq 1$, $\beta \equiv 0$ for $|x| \geq 2$, and $0 \leq \beta \leq 1$ everywhere. Let $q \in \R^2$ and define $\beta_q(x) = \beta(x - q)$. Hence, $\beta_q \equiv 1$ when $|x-q| \leq 1$ and $\beta_q \equiv 0$ when $|x-q| \geq 2$.

Define $\tilde{\beta_q} = \beta(\frac{x - q}{2})$ and $\tilde{\beta} = \beta(\frac{x}{2}) = \tilde{\beta_0}$. Hence, $\tilde{\beta_q} \equiv 1$ when $|x-q| \leq 2$ and $\tilde{\beta_q} \equiv 0$ when $|x-q| \geq 4$. Note that $\tilde{\beta_q} \equiv 1$ in $\supp \beta_q$.

Note that any of the above functions in $C_0^\infty(\R^2)$ can be considered functions in $\R^2 \times \R^2$ as well. In this case, they would be $C^\infty(\R^2 \times \R^2)$ but would no longer be compactly supported and would not depend on $\xi$.

Then there exists a finite collection of functions $\{ \chi_j \}_{j = 1}^n$ with the follow properties. 

\begin{enumerate}[a)] 
\item $\chi_1 \in C^\infty(\R^2 \times \R^2)$
\item $\chi_j \in C_0^\infty(\R^2 \times \R^2)$ for $j = 2, ..., n$
\item $\supp \chi_j \subseteq \supp \tilde{\beta}$ for all $j$
\item $0 \leq \chi_j \leq 1$ for all $j$
\item $\sum_{j = 1}^n \chi_j = \tilde{\beta}$
\item $|\tilde{p}_q| > \gamma > 0$ in $\supp \chi_1$ uniformly in $q$
\end{enumerate}

\begin{center}
\includegraphics[width = 0.7\textwidth]{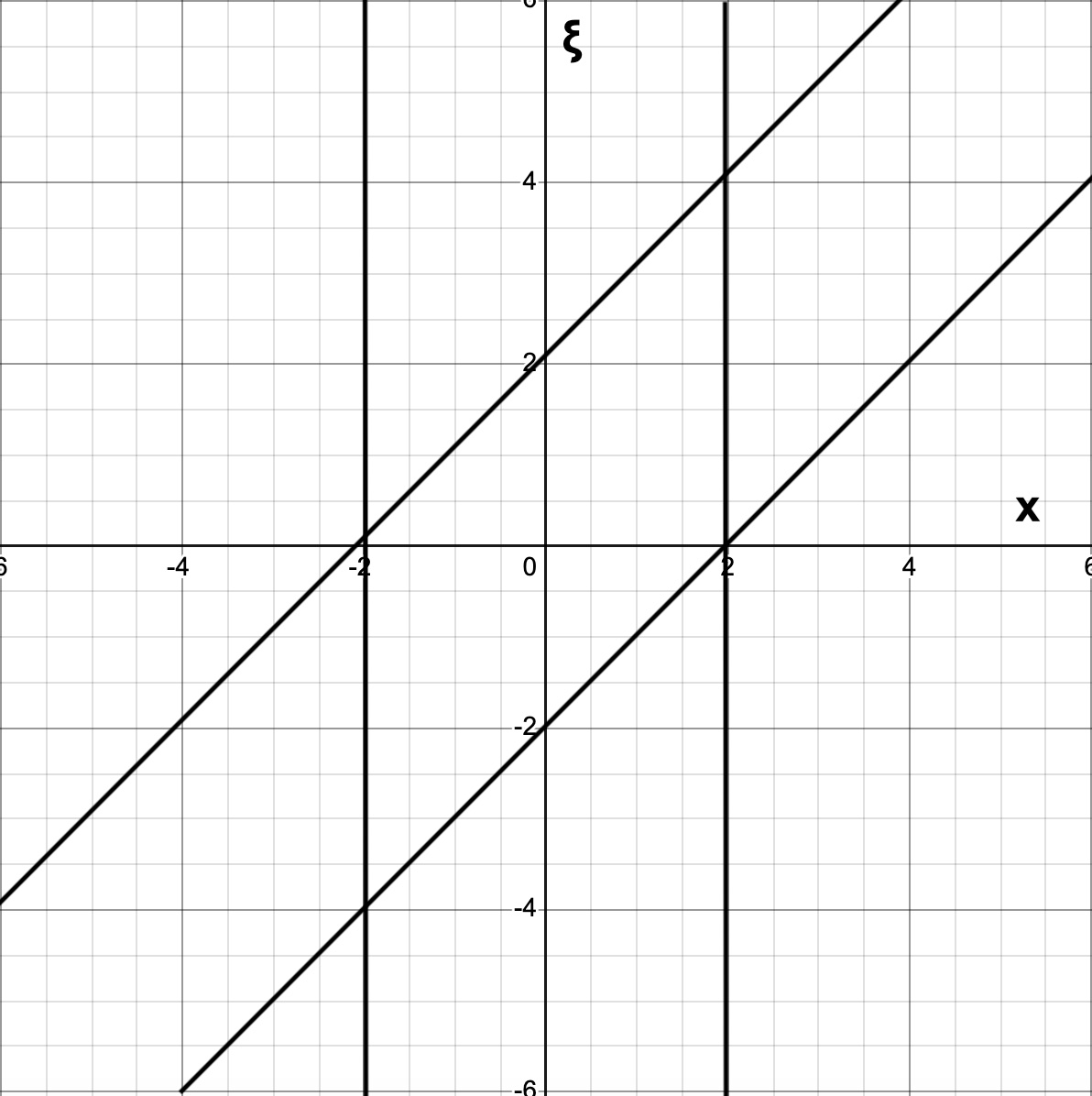}
\end{center}

We use the functions $\beta_q$ to localize in space and the functions $\chi_j$ to localize in frequency. The support of $\beta$ is contained in the vertical lines in the figure. The diagonal lines represent where $p(x, \xi) = 0$ for the specific case when $\phi(z) = |z|^2$, which can be thought of as the model case that we are generalizing. 

We begin with a few preparatory lemmas.

\begin{lemma} \thlabel{gextension_lemma}
    Let $\chi : \R^2 \times \R^2 \to \R$ such that $|x| < C$ in $\supp \chi$. Suppose there exists $0 < \gamma < 1$ such that $|\tilde{p}_q| > \gamma$ in $\supp \chi$. Then there exists $\tilde{\tilde{p}}_q \in C^\infty(\R^2 \times \R^2)$ such that

    \begin{enumerate}[a)] 
    \item $\tilde{\tilde{p}}_q \equiv \tilde{p}_q$ in $\supp \chi$
    \item $\tilde{\tilde{p}}_q \approx \langle \xi \rangle^2$, uniformly in $q$
    \item $\tilde{\tilde{p}}_q \in S_{0, 1/2}(\langle \xi \rangle^2)$, uniformly in $q$
    \item $|\tilde{\tilde{p}}_q| > \gamma/2$ everywhere
    \end{enumerate}

    \begin{proof}
        Without loss of generality, we assume that $\tilde{p}_q > \gamma > 0$ in $\supp \chi$. 

        Let $\tilde{\chi} \equiv 1$ in $\supp \chi$ and $\tilde{\chi} \equiv 0$ when $\text{dist}((x, \xi), \supp \chi) > \epsilon$ for some $\epsilon > 0$ small enough that $\tilde{p}_q > \gamma/2$ in $\supp \tilde{\chi}$. Note that this implies that $|x| < C + \epsilon$ in $\supp \tilde{\chi}$. Now put

        \[ \tilde{\tilde{p}}_q := \tilde{p}_q \tilde{\chi} + (1-\tilde{\chi}) \frac{\langle \xi \rangle^2}{4}\]

            \begin{enumerate}[a)] 
    \item That $\tilde{\tilde{p}}_q$ as so defined is identical to $\tilde{p}_q$ in $\supp \chi$ is obvious.
    
    \item  This is clear in the case that $\tilde{\chi} \equiv 1$ and $\tilde{\chi} \equiv 0$, and the intermediate case follows by elementary calculations.

    \item This follows from the Generalized Product Rule and the fact that $\langle \xi \rangle^2 \in S(\langle \xi \rangle^2)$ everywhere and $\tilde{p}_q \in S_{0,1/2}(\langle \xi \rangle^2)$ in $\supp \tilde{\chi}$, uniformly in $q$. 
    
    \item This follows from the fact that $\langle \xi \rangle^2 \geq 1 > \gamma$ everywhere and $p > \gamma/2$ in $\supp \tilde{\chi}$.
    \end{enumerate}

    \end{proof}
    
\end{lemma}

\begin{lemma} With $u_h, \beta_q,$ and $\chi_{j}$ defined as above, we have the following estimates.

\begin{equation} \label{ul3}
    ||u_h||_{L^\infty} \leq \sum_{j = 1}^n \sup_{q \in \Z^2} ||\chi_{j}^\mathrm{w}(x, hD) T_q \beta u_{h, -q}||_{L^\infty}
\end{equation}

and

\begin{equation} \label{ul3_L6}
    ||u_h||_{L^6} \leq \sum_{j = 1}^n \big( \sum_{q \in \Z^2} ||\chi_j^\mathrm{w}(x, hD) T_q \beta u_{h, -q}||_{L^6}^2 \big)^{1/2}
\end{equation}

\begin{proof}

We have

\begin{equation} \label{ul1}
    ||u_h||_{L^\infty(\R^2)} \leq \sup_{q \in \Z^2} ||\beta_q u_h||_{L^\infty(\R^2)} = \sup_{q \in \Z^2} ||T_q \beta u_{h, -q}||_{L^\infty(\R^2)} = \sup_{q \in \Z^2} || \tilde{\beta} T_q \beta u_{h, -q} ||_{L^\infty(\R^2)}
\end{equation}

We have used that $\tilde{\beta} \equiv 1$ in $\supp \beta$. Furthermore, because $\tilde{\beta}$ does not depend on $\xi$, its Weyl quantization is simply multiplication by the function $\tilde{\beta}$. Thus, $\tilde{\beta} T_q \beta u_{h, -q} = \tilde{\beta}^\mathrm{w} T_q \beta u_{h, -q}$. Furthermore, because $\tilde{\beta} = \sum_{j = 1}^n \chi_{j}$, by the linearity of the Weyl quantization,

\begin{equation} \label{ul2}
    \tilde{\beta} T_q \beta u_{h, -q} = \tilde{\beta}^\mathrm{w}(x, hD) \beta u_{h, -q} = (\sum_{j = 1}^n \chi_{j})^\mathrm{w}(x, hD) T_q \beta u_{h, -q} = \sum_{j = 1}^n \chi_{j}^\mathrm{w}(x, hD) T_q \beta u_{h, -q}
\end{equation}

Thus, by \eqref{ul1} and \eqref{ul2}, \eqref{ul3} is proved.

Likewise, we have

\begin{equation} \label{ul2_L6}
    ||u_h||_{L^6(\R^2)}^6 \lesssim \sum_{q \in \Z^2} ||\beta_q u_h||_{L^6(\R^2)}^6 = \sum_{q \in \Z^2} ||T_q \beta u_{h, -q}||_{L^6(\R^2)}^6 \lesssim \sum_{j = 1}^n \sum_{q \in \Z^2} ||\chi_j(x, hD)^\mathrm{w} T_q \beta u_{h, -q}||_{L^6(\R^2)}^6
\end{equation}

This implies \eqref{ul3_L6}, by the $l^1 \hookrightarrow l^3$ inequality.

\end{proof}

\end{lemma}

\begin{lemma} \thlabel{ABh_lemma}

For $u_h$ such that $P u_h = 0$, we have

\begin{equation*}
    ||Au_h||_{L^2(\R^2)}, ||Bu_h||_{L^2(\R^2)} \leq (\frac{h}{4} ||\Delta \phi||_{L^\infty(\R^2)} + 1)^{1/2}||u_h||_{L^2(\R^2)}
\end{equation*}

\begin{proof}
    We calculate

\begin{align*}
    ||Au_h||_{L^2(\R^2)}^2, ||Bu_h||_{L^2(\R^2)}^2 &\leq ||Au_h||_{L^2(\R^2)}^2 + ||Bu_h||_{L^2(\R^2)}^2 = \langle Au_h , Au_h \rangle + \langle Bu_h, Bu_h \rangle \\ &= \langle (A^2 + B^2) u_h, u_h \rangle = \langle (\frac{h^2}{4} \Delta \phi_h +1) u_h , u_h \rangle \leq ||\frac{h^2}{4} \Delta \phi_h + 1 ||_{L^\infty(\R^2)} ||u_h||_{L^2(\R^2)}^2 \\ &\leq (\frac{h}{4} ||\Delta \phi||_{L^\infty(\R^2)} + 1) ||u_h||_{L^2(\R^2)}^2
\end{align*}

where we have used that $A$ and $B$ are both self-adjoint and that $(A^2 + B^2)u_h = (\frac{h^2}{4} \Delta \phi_h +1)u_h$ because $Pu_h = 0$. 
    
\end{proof}

\end{lemma}

\begin{lemma} \thlabel{glcutoffefunc}

For $\beta \in C_0^\infty(\R^2)$, we have 

\begin{equation*}
    \sup_{q \in \R^2} ||P \beta_q u_h||_{L^2(\R^2)} \lesssim h||u_h||_{L^2(\R^2)}
\end{equation*}

and

\begin{equation*}
    \left( \sum_{q \in \Z^2} ||P \beta_q u_h||_{L^2(\R^2)}^2 \right)^{1/2} \lesssim h||u_h||_{L^2(\R^2)}
\end{equation*}

\begin{proof}


We can calculate directly, using the fact that $Pu_h = 0$:

\begin{align*}
    P \beta_q u_h = h \left( \frac{h}{4} (D_1^2 \beta_q + D_2^2 \beta_q) u_h + D_1 \beta_q (\frac{h}{2} D_1 - \frac{h}{2} \partial_2 \phi_h)u_h + D_2 \beta_q (\frac{h}{2} D_2 + \frac{h}{2} \partial_1 \phi_h) u_h \right) 
\end{align*}

Hence, we have

\[ \sup_{q \in \R^2} ||P \beta_q u_h||_{L^2(\R^2)} \leq h\left(\frac{h}{4}||\Delta \beta||_{L^\infty(\R^2)} ||u_h||_{L^2(\R^2)} + ||D_1 \beta||_{L^\infty(\R^2)} ||A u_h||_{L^2(\R^2)} + ||D_2 \beta||_{L^\infty(\R^2)} ||B u_h||_{L^2(\R^2)}\right)\]

and 

\begin{multline*}
    \left( \sum_{q \in \Z^2} ||P \beta_q u_h||_{L^2(\R^2)}^2 \right)^{1/2} \\ \lesssim h \left( \frac{h}{4} ||\sum_{q \in \Z^2} |\Delta \beta_q|^2||_{L^\infty}^{1/2}||u_h||_{L^2(\R^2)} + ||\sum_{q \in \Z^2} |D_1 \beta_q|^2||_{L^\infty}^{1/2}||Au_h||_{L^2(\R^2)} + ||\sum_{q \in \Z^2} |D_2 \beta_q|^2||_{L^\infty}^{1/2}||Bu_h||_{L^2(\R^2)}\right)
\end{multline*}

Hence, by \thref{ABh_lemma}, we arrive at

\[\sup_{q \in \R^2} ||P \beta_q u_h||_{L^2(\R^2)} \leq h \left(\frac{h}{4} ||\Delta \beta||_{L^\infty(\R^2)} + (\frac{h}{4} ||\Delta \phi||_{L^\infty(\R^2)} + 1)^{1/2} \left(||D_1 \beta||_{L^\infty(\R^2)} + ||D_2 \beta||_{L^\infty(\R^2)}\right)\right) ||u_h||_{L^2(\R^2)}\]

and

\begin{multline*}
    \left( \sum_{q \in \Z^2} ||P \beta_q u_h||_{L^2(\R^2)}^2 \right)^{1/2} \\ \lesssim h \left( \frac{h}{4} ||\sum_{q \in \Z^2} |\Delta \beta_q|^2||_{L^\infty}^{1/2} + (\frac{h}{4} ||\Delta \phi||_{L^\infty(\R^2)} + 1)^{1/2}(||\sum_{q \in \Z^2} |D_1 \beta_q|^2||_{L^\infty}^{1/2} + ||\sum_{q \in \Z^2} |D_2 \beta_q|^2||_{L^\infty}^{1/2})\right) ||u_h||_{L^2(\R^2}
\end{multline*}

Note that the $L^\infty$ norms in the above expression are finite due to the finite overlap of the functions $\beta_q$. 

\end{proof}

\end{lemma}

\begin{lemma} \thlabel{L_infty_to_L2}
    \begin{equation}
        ||u||_{L^\infty(\R^2)} \lesssim h^{-1} ||\langle hD \rangle^2 u||_{L^2(\R^2)}
    \end{equation}

    and 

    \begin{equation} \thlabel{L6_to_L2}
        ||u||_{L^6(\R^2)} \lesssim h^{-1 + 1/3} ||\langle hD \rangle^2 u||_{L^2(\R^2)}
    \end{equation}

    \begin{proof}
        \begin{align*}
            ||u||_{L^\infty} &\lesssim h^{-2} ||\mathcal{F}_h u||_{L^1} \\ &= h^{-2} ||\langle \cdot \rangle^{-2} \langle \cdot \rangle^2 \mathcal{F}_h u||_{L^1} \\ &\leq h^{-2} ||\langle \cdot \rangle^{-2}||_{L^2} ||\langle \cdot \rangle^2 \mathcal{F}_h u||_{L^2} \\ &\lesssim h^{-1} || \mathcal{F}_h \langle \xi \rangle^2 \mathcal{F}_h u||_{L^2} \\ &= h^{-1} ||\langle hD \rangle^2 u||_{L^2}
        \end{align*}

        The proof of \eqref{L6_to_L2} is essentially the same except we interpolate between

        \begin{equation*}
            ||u||_{L^\infty} \lesssim h^{-2} ||\mathcal{F}_h u||_{L^1}
        \end{equation*}

        and

        \begin{equation*}
            ||u||_{L^2} \lesssim h^{-1} ||\mathcal{F}_h u||_{L^2}
        \end{equation*}

        to find that

        \begin{equation*}
            ||u||_{L^6} \lesssim h^{-2 + 1/3} ||\mathcal{F}_h u||_{L^{6/5}}
        \end{equation*}

        and use a generalized Hölder inequality.

    \end{proof}
\end{lemma}

Now, we split the theorem into two cases. In the first case, we have $\tilde{p}_q$ may be equal to 0 at some points in the support of $\chi_j$. Note that we can choose our partition in such a way that $|\tilde{p}_q|$ is small in the support of these $\chi_j$. After this case, we deal with the case when $\tilde{p}_q$ is bounded away from 0 in the support of $\chi_1$. This second case is essentially an elliptic estimate.

\begin{theorem}
    We have the following estimates uniformly in $q$ for $h$ sufficiently small when $\tilde{p}_q = 0$ at some points in $\supp \chi_{j}$

    \begin{equation}
        ||\chi_{j}^\mathrm{w}(x, hD) T_q \beta u_{h, -q}||_{L^\infty(\R^2)} \lesssim h^{-1/2} \left(||\beta_q u_h||_{L^2(\R^2)} + \frac{1}{h}||P \beta_q u_h||_{L^2(\R^2)} \right)
    \end{equation}

    Likewise, we have

    \begin{equation} \label{non-elliptic_L6}
        ||\chi_{j}^\mathrm{w}(x, hD)  T_q \beta u_{h, -q}||_{L^6(\R^2)} \lesssim h^{-1/6} \left( ||\beta_q u_h||_{L^2(\R^2)} + \frac{1}{h} ||P \beta_q u_h||_{L^2(\R^2)} \right)
    \end{equation}

    \begin{proof}
        Recall that $\tilde{p}_q(x, \xi) = \tilde{a}_q^2 + \tilde{b}_q^2 - (\frac{h^2}{4}(\Delta \phi_h)(x)+1)$. Without loss of generality, we can suppose that, in the support of $\chi_{j}$, $\tilde{a}_q > \frac{1}{\sqrt{2}} - \epsilon$ and $|\tilde{b}_q| < \frac{1}{\sqrt{2}} + \epsilon$, which corresponds to the top rectangle in the figure. The $\chi_j$ corresponding to the other rectangles in the figure can be handled analogously.
        
\begin{center}
\includegraphics[width = 0.7\textwidth]{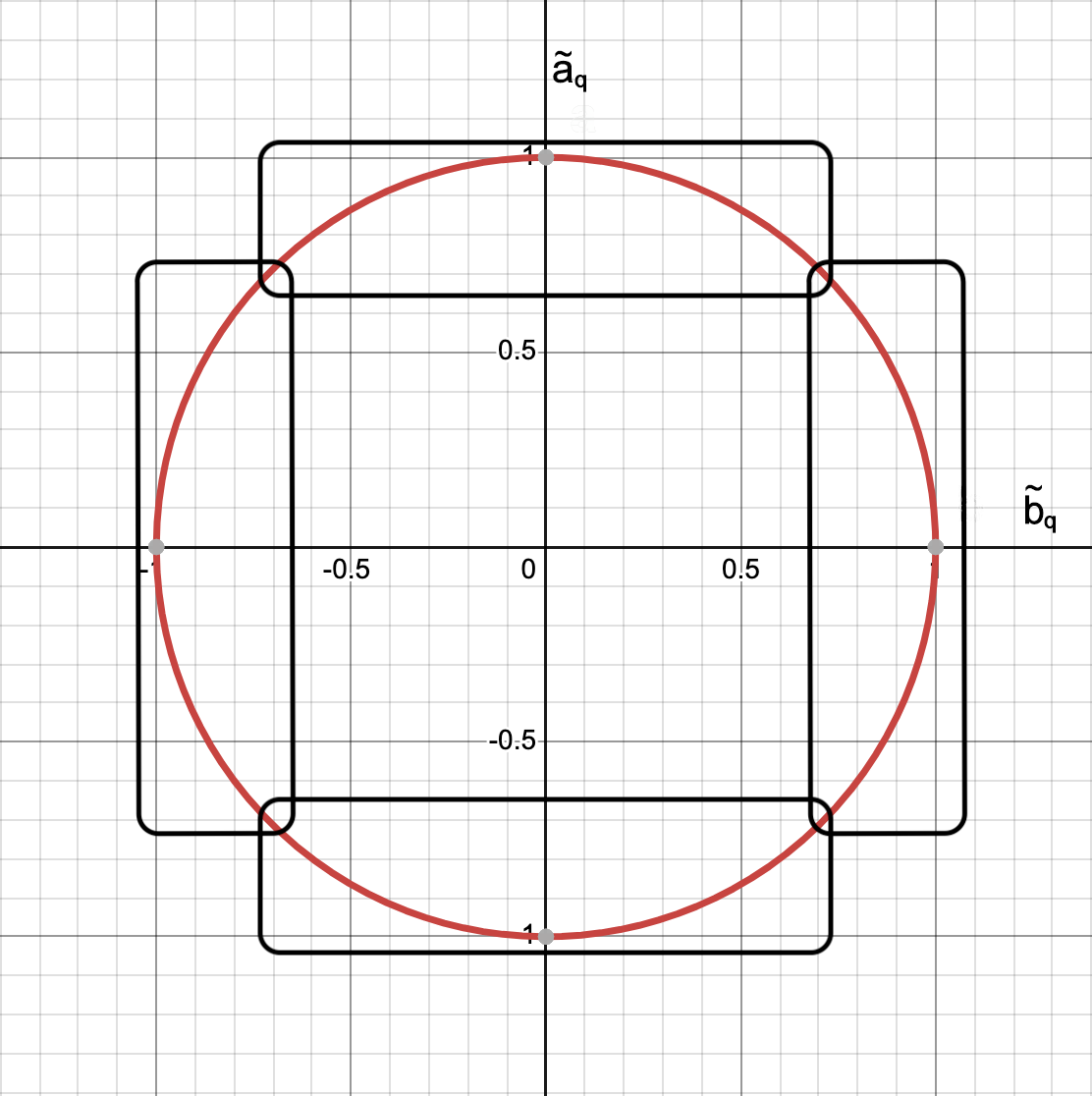}
\end{center}
        
        Now, let $\tilde{\psi}_{j}, \tilde{\tilde{\psi}}_{j} \in C_0^\infty(\R^2 \times \R^2)$ such that $\tilde{\psi}_{j} \equiv 1$ in the support of $\chi_{j}$ and $\tilde{\tilde{\psi}}_{j} \equiv 1$ in the support of $\tilde{\psi}_{j}$. We can also suppose that $\tilde{a}_q > \frac{1}{\sqrt{2}} - \epsilon$ and $|\tilde{b}_q| < \frac{1}{\sqrt{2}} + \epsilon$ in the supports of $\tilde{\psi}_{j}$ and $\tilde{\tilde{\psi}}_{j}$.

        Now, define

        \begin{align*}
            s(x, \xi) &:= \frac{\xi_1}{2} + \left(-\frac{h}{2} \partial_2 \phi_h(x+q) + \frac{h}{2} \partial_2 \phi_h(q)\right)\tilde{\psi}_j - \sqrt{(\frac{h^2}{4}(\Delta \phi_h)(x)+1 - \tilde{b}_q^2) \tilde{\psi}_{j}^2} \\
            \tilde{s}(x, \xi) &:= \tilde{a}_q\tilde{\psi}_{j} - \sqrt{(\frac{h^2}{4}(\Delta \phi_h)(x)+1 - \tilde{b}_q^2) \tilde{\psi}_{j}^2} \\
        \end{align*}

        We likewise define

        \begin{equation*}
            e(x, \xi) := \tilde{a}_q \tilde{\tilde{\psi}}_{j} + \sqrt{(\frac{h^2}{4}(\Delta \phi_h)(x)+1 - \tilde{b}_q^2) \tilde{\tilde{\psi}}_{j}^2}
        \end{equation*}

        in the support of $\tilde{\psi}_j$ and extend $e$ arbitrarily so that $e \in S_{0, 1/2}(1)$ and $e$ is bounded away from 0. Thus, $e$ is an elliptic operator:

        \begin{equation*}
            ||e^\mathrm{w}(x, hD) u||_{L^2} \gtrsim ||u||_{L^2}
        \end{equation*}

        for $u \in \mathcal{S}$ (see Theorem 4.29 in \cite{Zworski}). 
        
        By a Bernstein inequality and an energy estimate (\thref{Bernstein} and \thref{energy_estimate}, see the appendix), we have

        \begin{align*}
            ||\chi_j^\mathrm{w} T_q \beta u_{h, -q}||_{L^\infty(\R^2)} &\lesssim h^{-1/2} ||\chi_j^\mathrm{w} T_q \beta u_{h, -q}||_{L^\infty(\R) L^2(\R)} \\ &\lesssim h^{-3/2}||s^\mathrm{w}\chi_j^\mathrm{w} T_q \beta u_{h, -q}||_{L^1(\R) L^2(\R^)}\\ &\lesssim h^{-3/2}||s^\mathrm{w}\chi_j^\mathrm{w} T_q \beta u_{h, -q}||_{L^2(\R^2)}
        \end{align*}

        where we have used the compact support in $x_1$ to set an initial condition of 0 and dominate the $L^1(\R)$ norm by the $L^2(\R)$ norm. 

        Likewise, if we let
        
        \begin{equation*}
            f(x_1, x_2) = s^\mathrm{w} \chi_j^\mathrm{w} T_q \beta u_{h, -q}
        \end{equation*}
        
        we have, for some $t_1, t_2$ and solution operator $U$,

        \begin{equation*}
            ||\chi_j^\mathrm{w} T_q \beta u_{h, -q}||_{L^6(\R^2)} \lesssim ||\frac{1}{h} \int_{t_1}^{t_2} |U(t,s) f(s, x_2)| ds ||_{L^6(\R^2)} \lesssim h^{-1/6} \left( \frac{1}{h} ||s^\mathrm{w} \chi_j^\mathrm{w} T_q \beta u_{h, -q}||_{L^2(\R^2)} \right)
        \end{equation*}

        This follows from Duhamel's Formula and a microlocal version of the Strichartz estimate (\thref{strichartz}, see the appendix) applied to the Schrödinger operator $s^\mathrm{w}(x, hD)$.

        We note that we can show that the assumption of \thref{strichartz} is satisfied by direct calculation of the second derivative of

        \begin{equation*}
            (-\frac{h}{2} \partial_2 \phi_h(x+q) + \frac{h}{2} \partial_2 \phi_h(q))\tilde{\psi}_j - \sqrt{(\frac{h^2}{4}(\Delta \phi_h)(x)+1 - \tilde{b}_q^2) \tilde{\psi}_{j}^2}
        \end{equation*}

        with respect to $\xi_2$ in the support of $\chi_j$. By direct calculation, we see that this equals

        \begin{equation*}
            \frac{1}{4} \left[ \frac{\frac{h^2}{4} \Delta \phi_h + 1}{(\frac{h^2}{4} \Delta \phi_h + 1 - \tilde{b}_q^2)^{3/2}} \right]
        \end{equation*}

        which is clearly non-degenerate in the support of $\chi_j$.
        
        Because $s - \tilde{s} = \frac{\xi_1}{2} (1 - \tilde{\psi}_j)$, $\tilde{\psi}_j \equiv 1$ in $\supp \chi_j$, and $\chi_j \in C_0^\infty$ we have that

        \begin{equation*}
            ||(s - \tilde{s})^\mathrm{w} \chi_j^\mathrm{w} T_q \beta u_{h, -q}||_{L^2(\R^2)} \lesssim h^\infty||\beta_q u_h||_{L^2(\R^2)}
        \end{equation*}

        Furthermore, because $e$ is elliptic, we have

        \begin{equation*}
            ||\tilde{s}^\mathrm{w} \chi_j^\mathrm{w} T_q \beta u_{h, -q}||_{L^2(\R^2)} \lesssim ||e^\mathrm{w} \tilde{s}^\mathrm{w} \chi_j^\mathrm{w} T_q \beta u_{h, -q}||_{L^2(\R^2)}
        \end{equation*}

        Now, note that

        \begin{equation}
            e(x, \xi) s(x, \xi) = \tilde{p}_q(x, \xi) \tilde{\psi}_{j} \tilde{\tilde{\psi}}_{j} = \tilde{p}_q(x, \xi) \tilde{\psi}_{j} 
        \end{equation}

        Thus,

        \begin{align*}
            ||e^\mathrm{w} \tilde{s}^\mathrm{w} \chi_j^\mathrm{w} T_q \beta u_{h, -q}||_{L^2(\R^2)} &\leq ||(e^\mathrm{w} \tilde{s}^\mathrm{w} - (e\tilde{s})^\mathrm{w}) \chi_j^\mathrm{w} T_q \beta u_{h, -q}||_{L^2(\R^2)} \\ &+ ||((\tilde{p_q} \tilde{\psi}_j)^\mathrm{w} - \tilde{p}_q^\mathrm{w} \tilde{\psi}_j^\mathrm{w})\chi_j^\mathrm{w} T_q \beta u_{h, -q}||_{L^2(\R^2)} \\ &+ ||\tilde{p}_q^\mathrm{w} (\tilde{\psi}_j^\mathrm{w} 
 - 1)\chi_j^\mathrm{w} T_q \beta u_{h,-q}||_{L^2(\R^2)} \\ &+ ||\tilde{p}_q^\mathrm{w} \chi_j^\mathrm{w} T_q \beta u_{h, -q}||_{L^2(\R^2)}
        \end{align*}

        We consider each of these summands in turn. Because

        \begin{equation}
            e^\mathrm{w}\tilde{s}^\mathrm{w} - (e\tilde{s})^\mathrm{w} = \frac{h}{2i} \{e, \tilde{s}\}^\mathrm{w} + O(h^{3/2})
        \end{equation}

        and

        \begin{equation}
            \tilde{p}_q^\mathrm{w}\tilde{\psi}_j^\mathrm{w} - (\tilde{p}_q\tilde{\psi}_j)^\mathrm{w} = \frac{h}{2i} \{\tilde{p}_q, \tilde{\psi}_j\}^\mathrm{w} + O(h^{3/2})
        \end{equation}

        we have

        \begin{equation}
            ||(e^\mathrm{w} \tilde{s}^\mathrm{w} - (e \tilde{s})^\mathrm{w})\chi_j^\mathrm{w} T_q \beta u_{h, -q}||_{L^2(\R^2)} \lesssim h||\beta_q u_h||_{L^2(\R^2)}
        \end{equation}

        and 

        \begin{equation}
            ||(\tilde{p}_q^\mathrm{w} \tilde{\psi_j}^\mathrm{w} - (\tilde{p}_q \tilde{\psi}_j)^\mathrm{w})\chi_j^\mathrm{w} T_q \beta u_{h, -q}||_{L^2(\R^2)} \lesssim h||\beta_q u_h||_{L^2(\R^2)}
        \end{equation}

        Likewise, because $\tilde{\psi}_j \equiv 1$ in $\supp \chi_j$, we have

        \begin{equation}
            ||\tilde{p}_q^\mathrm{w} (\tilde{\psi}_j^\mathrm{w} - 1)\chi_j^\mathrm{w} T_q \beta u_{h, -q}||_{L^2(\R^2)} \lesssim h^\infty ||\beta_q u_h||_{L^2}
        \end{equation}

    Finally,

    \begin{equation*}
        ||\tilde{p}_q^\mathrm{w} \chi_j^\mathrm{w} T_q \beta u_{h, -q}||_{L^2(\R^2)} \lesssim ||[\tilde{p}_q^\mathrm{w}, \chi_j^\mathrm{w}]T_q \beta u_{h,-q}||_{L^2(\R^2)} + ||\tilde{p}_q^\mathrm{w} T_q \beta u_{h, -q}||_{L^2(\R^2)}
    \end{equation*}

    Now, because

    \begin{equation*}
        [\tilde{p}_q^\mathrm{w}, \chi_j^\mathrm{w}] = \frac{h}{i} \{\tilde{p}_q, \chi_j \}^\mathrm{w} + O(h^{3/2})
    \end{equation*}

    we have

    \begin{equation*}
        ||[\tilde{p}_q^\mathrm{w}, \chi_j^\mathrm{w}]T_q \beta u_{h,-q}||_{L^2(\R^2)} \lesssim h||\beta_q u_h||_{L^2(\R^2)}
    \end{equation*}

    and

    \begin{equation*}
        ||\tilde{p}_q^\mathrm{w}T_q \beta u_{h, -q}||_{L^2(\R^2)} = ||T_q^{-1} \tilde{p}_q^\mathrm{w} T_q \beta u_{h, -q}||_{L^2(\R^2)} = ||p^\mathrm{w} \beta_q u_h||_{L^2(\R^2)}
    \end{equation*}

%
%
%

    \end{proof}

\end{theorem}

\begin{theorem}
    With $u_h, \beta, \chi_{j},$ and $p$ defined as above, we have the following estimate uniformly in $q$ for $h$ sufficiently small, provided that there exists $\gamma > 0$ such that $|\tilde{p}_q| > \gamma$ in $\supp \chi_{1}$.

    \begin{equation}
        ||\chi_{1}^\mathrm{w} T_q \beta u_{h, -q}||_{L^\infty(\R^2)} \lesssim h^{-1/2} \left(||\beta_q u_{h}||_{L^2(\R^2)} + \frac{1}{h} ||P \beta_q u_h||_{L^2(\R^2)} \right)
    \end{equation}

    Likewise, we have

    \begin{equation}
        ||\chi_{1}^\mathrm{w} T_q \beta u_{h, -q}||_{L^6(\R^2)} \lesssim h^{-1/6} \left(||\beta_q u_h||_{L^2(\R^2)} + \frac{1}{h} ||P \beta_q u_h||_{L^2(\R^2)} \right)
    \end{equation}

    \begin{proof}

        By \thref{L_infty_to_L2}, it suffices to show that

        \begin{equation*}
            ||\langle hD \rangle^2 \chi_1^\mathrm{w} T_q \beta u_{h, -q}||_{L^2(\R^2)} \lesssim h^{1/2} \left( ||\beta_q u_h||_{L^2(\R^2)} + \frac{1}{h} ||P \beta_q u_h||_{L^2(\R^2)} \right)
        \end{equation*}

        By \thref{gextension_lemma}, there exists $\tilde{\tilde{p}}_q \in S_{0,1/2}(\langle \xi \rangle^2)$ such that $\tilde{\tilde{p}}_q \equiv \tilde{p}_q$ in $\supp \chi_{1}, \tilde{\tilde{p}}_q \approx \langle \xi \rangle^2$ and $|\tilde{\tilde{p}}_q| > \gamma/2$ everywhere. Hence, there exists $g_q \in S_{0,1/2}(1/\langle \xi \rangle^2)$, uniformly in $q$ such that $g_q^\mathrm{w} \tilde{\tilde{p}}_q^\mathrm{w} = I + h^N r_{N, q}^\mathrm{w}$ where $r_{N, q} \in S_{1/2}$, uniformly in $q$, for any $N \in \N$. 

        Thus,

        \begin{align*}
            ||\langle hD \rangle^2 \chi_1^\mathrm{w} T_q \beta u_{h, -q}||_{L^2(\R^2)} &= ||\langle hD \rangle^2 (g_q^\mathrm{w} \tilde{\tilde{p}}_q^\mathrm{w} - h^N r_{N, q}^\mathrm{w}) \chi_1^\mathrm{w} T_q \beta u_{h, -q}||_{L^2(\R^2)} \\ &\lesssim ||\tilde{\tilde{p}}_q^\mathrm{w} \chi_1^\mathrm{w} T_q \beta u_{h, - q}||_{L^2(\R^2)} + h^N ||\langle hD \rangle^2 r_{N, q}^\mathrm{w} \chi_1^\mathrm{w} T_q \beta u_{h, -q}||_{L^2(\R^2)}
        \end{align*}

%
%
%
%
%

        
        Using the explicit construction of $\tilde{\tilde{p}}_q$, we calculate that

        \begin{align*}
            \tilde{\tilde{p}}_q - \tilde{p}_q &= \left(1 - \tilde{\chi}\right)\left(\frac{\langle \xi \rangle^2}{4} - \tilde{p}_q \right) \\ &= (1 - \tilde{\chi})\left(\frac{5}{4} + \frac{h^2}{4} \Delta \phi_h - \left(\frac{h}{2} \partial_2 \phi_h(x+q) - \frac{h}{2} \partial_2 \phi_h(q)\right)^2 - \left(\frac{h}{2} \partial_1 \phi_h(x+q) - \frac{h}{2} \partial_1 \phi_h(q)\right)^2 \right. \\ & \left. + \xi_1\left(\frac{h}{2} \partial_2 \phi_h(x+q) - \frac{h}{2} \partial_2 \phi_h(q)\right) - \xi_2\left(\frac{h}{2} \partial_1 \phi_h(x+q) - \frac{h}{2} \partial_1 \phi_h(q)\right) \right)
        \end{align*}

        Let 

        \begin{align*}
            v_1(x, \xi) &:= (1 - \tilde{\chi})\left(\frac{h}{2} \partial_2 \phi_h(x+q) - \frac{h}{2} \partial_2 \phi_h(q)\right)^2 \\
            v_2(x, \xi) &:= (1 - \tilde{\chi})\left(\frac{h}{2} \partial_1 \phi_h(x+q) - \frac{h}{2} \partial_1 \phi_h(q)\right)^2 \\
            v_3(x, \xi) &:= (1 - \tilde{\chi})\xi_1 \left(\frac{h}{2} \partial_2 \phi_h(x+q) - \frac{h}{2} \partial_2 \phi_h(q)\right) \\
            v_4(x, \xi) &:= (1 - \tilde{\chi})\xi_2 \left(\frac{h}{2} \partial_1 \phi_h(x+q) - \frac{h}{2} \partial_1 \phi_h(q)\right)
        \end{align*}

        Hence, 

        \begin{align*}
            ||(\tilde{\tilde{p}}_q - \tilde{p}_q)^\mathrm{w} \chi_1^\mathrm{w} T_q \beta u_{h, -q}||_{L^2(\R^2)} &\lesssim ||\left((1 - \tilde{\chi})\left(\frac{5}{4} + \frac{h^2}{4} \Delta \phi_h\right)\right)^\mathrm{w} \chi_1^\mathrm{w} T_q \beta u_{h, -q}||_{L^2(\R^2)} \\ &+ ||v_1^\mathrm{w} \chi_1^\mathrm{w} T_q \beta u_{h, -q}||_{L^2(\R^2)} +  ||v_2^\mathrm{w} \chi_1^\mathrm{w} T_q \beta u_{h, -q}||_{L^2(\R^2)} \\ &+  ||v_3^\mathrm{w} \chi_1^\mathrm{w} T_q \beta u_{h, -q}||_{L^2(\R^2)} + ||v_4^\mathrm{w}\chi_1^\mathrm{w}T_q \beta u_{h, -q}||_{L^2(\R^2)}
        \end{align*}

        Because $1 - \tilde{\chi} \equiv 0$ in $\supp \chi_1$ and $(1 - \tilde{\chi})(\frac{5}{4} + \frac{h^2}{4} \Delta \phi_h) \in S_{1/2}$, we have

        \begin{equation*}
             ||\left((1 - \tilde{\chi})\left(\frac{5}{4} + \frac{h^2}{4} \Delta \phi_h\right)\right)^\mathrm{w} \chi_1^\mathrm{w} T_q \beta u_{h, -q}||_{L^2(\R^2)} \lesssim h^\infty ||\beta_q u_h||_{L^2(\R^2)}
        \end{equation*}

        Likewise, because $v_1, v_2 \equiv 0$ in $\supp \chi_1$ and $v_1, v_2 \in S_{1/2}(\langle x \rangle^2)$, and because $\chi_1 \in S(\langle x \rangle^{-2})$, we have

        \begin{align*}
            ||v_1^\mathrm{w} \chi_1^\mathrm{w} T_q \beta u_{h, -q}||_{L^2(\R^2)}, ||v_2^\mathrm{w} \chi_1^\mathrm{w} T_q \beta u_{h, -q}||_{L^2(\R^2)} \lesssim h^\infty ||\beta_q u_h||_{L^2(\R^2)}
        \end{align*}

        We furthermore have

        \begin{align*}
            ||v_3^\mathrm{w} \chi_1^\mathrm{w} T_q \beta u_{h, -q}||_{L^2(\R^2)} &\lesssim ||\left(v_3^\mathrm{w} - hD_1 \left((1-\tilde{\chi})\left(\frac{h}{2} \partial_2 \phi_h(x+q) - \frac{h}{2} \partial_2 \phi_h(q)\right)\right)^\mathrm{w}\right) \chi_1^\mathrm{w} T_q \beta u_{h, -q}||_{L^2(\R^2)} \\ &+ ||hD_1 \left((1-\tilde{\chi})\left(\frac{h}{2} \partial_2 \phi_h(x+q) - \frac{h}{2} \partial_2 \phi_h(q)\right)\right)^\mathrm{w} \chi_1^\mathrm{w} T_q \beta u_{h, -q}||_{L^2(\R^2)}
        \end{align*}

        Then

        \begin{align*}
            v_3^\mathrm{w} - hD_1 \left((1-\tilde{\chi})\left(\frac{h}{2} \partial_2 \phi_h(x+q) - \frac{h}{2} \partial_2 \phi_h(q)\right)\right)^\mathrm{w} &= h\left\{ \xi_1, (1-\tilde{\chi})\left(\frac{h}{2} \partial_2 \phi_h(x+q) - \frac{h}{2} \partial_2 \phi_h(q)\right) \right\}^\mathrm{w} \\ &= h \left(\partial_{x_1} \left((1-\tilde{\chi})\left(\frac{h}{2} \partial_2 \phi_h(x+q) - \frac{h}{2} \partial_2 \phi_h(q)\right)\right)\right)^\mathrm{w}
        \end{align*}

        Note this follows from the asymptotic expansion of the symbol of $hD_1 ((1-\tilde{\chi})(\frac{h}{2} \partial_2 \phi_h(x+q) - \frac{h}{2} \partial_2 \phi_h(q)))^\mathrm{w}$, which in fact consists of only two terms because any second and higher derivatives of $\xi_1$ are 0. This symbol still has disjoint support from $\chi_1$ and now does not grow in $\xi$. Thus, we have

        \begin{align*}
            ||\left(v_3^\mathrm{w} - hD_1 \left((1-\tilde{\chi})\left(\frac{h}{2} \partial_2 \phi_h(x+q) - \frac{h}{2} \partial_2 \phi_h(q)\right)\right)^\mathrm{w}\right) \chi_1^\mathrm{w} T_q \beta u_{h, -q}||_{L^2(\R^2)} \lesssim h^\infty ||\beta_q u_h||_{L^2(\R^2)}
        \end{align*}

        Now, because $(1-\tilde{\chi})(\frac{h}{2} \partial_2 \phi_h(x+q) - \frac{h}{2} \partial_2 \phi_h(q)) \in S_{1/2}(\langle x \rangle)$ and $\chi_1 \in S(\langle x \rangle^{-1})$ and they have disjoint supports, we have, for any $N \in \N$, there exists $s_1 \in S_{1/2}$ such that

        \begin{align*}
             &h^{-1}||hD_1 \left((1-\tilde{\chi})\left(\frac{h}{2} \partial_2 \phi_h(x+q) - \frac{h}{2} \partial_2 \phi_h(q)\right)\right)^\mathrm{w} \chi_1^\mathrm{w} T_q \beta u_{h, -q}||_{L^2(\R^2)} = h^N ||hD_1 s_1^\mathrm{w} T_q \beta u_{h,-q}||_{L^2(\R^2)} \\ &\leq h^N ||hD_1 T_q \beta u_{h, -q}||_{L^2(\R^2)} + h^N ||[hD_1, s^\mathrm{w}] T_q \beta u_{h, -q}||_{L^2(\R^2)} \\ &\lesssim h^N||\tilde{A}_q T_q \beta u_{h, -q}||_{L^2(\R^2)} + h^N||\left(\frac{h}{2} \partial_2 \phi_h (x+q) - \frac{h}{2} \partial_2 \phi_h(q)\right) T_q \beta u_{h, -q}||_{L^2(\R^2)} + h^N ||\beta_q u_h||_{L^2(\R^2)}
        \end{align*}

        Note that we used \thref{commutator_lemma} above.

        Then

        \begin{align*}
            h^N ||\tilde{A}_q T_q \beta u_{h, -q}||_{L^2(\R^2)} = h^N ||T_{q}^{-1} \tilde{A}_q T_q \beta u_{h, -q}||_{L^2(\R^2)} = h^N ||A \beta_q u_{h}||_{L^2(\R^2)}
        \end{align*}

        and

        \begin{align*}
            h^N||\left(\frac{h}{2} \partial_2 \phi_h (x+q) - \frac{h}{2} \partial_2 \phi_h(q)\right) T_q \beta u_{h, -q}||_{L^2(\R^2)} \lesssim h^N ||\beta_q u_{h}||_{L^2(\R^2)}
        \end{align*}

        A similar estimate can be shown for the $v_4$ term. 
        
        Now,

        \begin{align*}
            ||\tilde{p}_q^\mathrm{w} \chi_1^\mathrm{w} T_q \beta u_{h, -q}||_{L^2(\R^2)} \leq  ||\tilde{p}_q^\mathrm{w} T_q \beta u_{h, -q}||_{L^2(\R^2)} + \sum_{j = 2}^n ||\tilde{p}_q^\mathrm{w} \chi_j^\mathrm{w} T_q \beta u_{h, -q}||_{L^2(\R^2)}
        \end{align*}

        because $\chi_1 = \tilde{\beta} - \sum_{j =2}^n \chi_j$ and $\tilde{\beta} \equiv 1$ in $\supp \beta$. Then

        \begin{equation*}
            ||\tilde{p}_q^\mathrm{w} T_q \beta u_{h, -q}||_{L^2(\R^2)} \lesssim ||p^\mathrm{w} \beta_q u_h||_{L^2(\R^2)} + h||\beta_q u_h||_{L^2(\R^2)}
        \end{equation*}

        and

        \begin{equation*}
            ||\tilde{p}_q^\mathrm{w} \chi_j^\mathrm{w} T_q \beta u_{h, -q}||_{L^2(\R^2)} \lesssim ||\tilde{p}_q^\mathrm{w} T_q \beta u_{h, -q}||_{L^2(\R^2)} + h ||\beta_q u_h||_{L^2(\R^2)}
        \end{equation*}

        for $j = 2,..., n$ because $\chi_j \in C_0^\infty(\R^2 \times \R^2)$.

        Finally, we turn to the term $h^N ||\langle hD \rangle^2 (r_{N, q} \# \chi_j)^\mathrm{w} T_q \beta u_{h, -q}||_{L^2(\R^2)}$. Note that we don't have to be too careful with the $h$'s because we have as many as we want. 

        \begin{align*}
            ||\langle hD \rangle^2 (r_{N, q} \# \chi_j)^\mathrm{w} T_q \beta u_{h, -q}||_{L^2(\R^2)} \lesssim ||\beta_q u_h||_{L^2(\R^2)} + ||\left((hD_1)^2 + (hD_2)^2\right) (r_{N, q} \# \chi_j)^\mathrm{w} T_q \beta u_{h, -q}||_{L^2(\R^2)}
        \end{align*}

        Furthermore,

        \begin{align*}
            ||\left((hD_1)^2 + (hD_2)^2\right) (r_{N, q} \# \chi_j)^\mathrm{w} T_q \beta u_{h, -q}||_{L^2(\R^2)} &\lesssim ||\left((hD_1)^2 + (hD_2)^2\right) T_q \beta u_{h, -q}||_{L^2(\R^2)} \\ &+ ||[(hD_1)^2 + (hD_2)^2, (r_{N, q} \# \chi_j)^\mathrm{w} ] T_q \beta u_{h, -q}||_{L^2(\R^2)}
        \end{align*}

        We handle the first term by adding and subtracting terms so that we have

        \begin{align*}
            ||\left((hD_1)^2 + (hD_2)^2\right) T_q \beta u_{h, -q}||_{L^2(\R^2)} &\lesssim ||\tilde{p}_q^\mathrm{w} T_q \beta u_{h, -q}||_{L^2(\R^2)} \\ &+ ||\tilde{A}_q T_q \beta u_{h, -q}||_{L^2(\R^2)} + || \tilde{B}_q T_q \beta u_{h, -q}||_{L^2(\R^2)} + ||\beta_q u_{h}||_{L^2(\R^2)}
        \end{align*}

            Finally, we handle the second term, using the fact that
            
            \begin{align*}
                [(hD_1)^2 + (hD_2)^2, (r_{N, q} \# \chi_j)^\mathrm{w}] &= [hD_1, [hD_1, (r_{N, q} \# \chi_j)^\mathrm{w}]] + 2 [hD_1, (r_{N, q} \# \chi_j)^\mathrm{w}] hD_1 \\ &+ [hD_2, [hD_2, (r_{N, q} \# \chi_j)^\mathrm{w}]] + 2 [hD_2, (r_{N, q} \# \chi_j)^\mathrm{w}] hD_2  
            \end{align*}

    combined with \thref{commutator_lemma} and \thref{L2_bound_for_S}. Thus, we have

    \begin{align*}
        ||[(hD_1)^2 + (hD_2)^2, (r_{N, q} \# \chi_j)^\mathrm{w} ] T_q \beta u_{h, -q}||_{L^2(\R^2)} &\lesssim ||\beta_q u_h||_{L^2(\R^2)} + ||A \beta_q u_h||_{L^2(\R^2)} + ||B \beta_q u_h||_{L^2(\R^2)}
    \end{align*}
    \end{proof}
\end{theorem}

\section{Small Eigenvalue Case}

In this section, we turn to the small eigenvalue case. As the $L^6$ estimate follows directly from the $L^\infty$ estimate, we concern ourselves in this section only with the latter. 

Define

\begin{align*}
    A &:= \frac{D_1}{2} - \frac{\partial_2 \phi(x)}{2} \\
    B &:= \frac{D_2}{2} + \frac{\partial_1 \phi(x)}{2} \\
    P &:= (\frac{D_1}{2} - \frac{\partial_2 \phi(x)}{2})^2 + (\frac{D_2}{2} + \frac{\partial_1 \phi(x)}{2})^2 - \frac{\Delta \phi(x)}{4} - \lambda^2 = A^2 + B^2 - \frac{\Delta \phi(x)}{4} - \lambda^2 \\
    A_q &:= \frac{D_1}{2} - \frac{\partial_2 \phi(x - q)}{2} \\
    B_q &:= \frac{D_2}{2} + \frac{\partial_1 \phi(x - q)}{2} \\
    P_q &:= (\frac{D_1}{2} - \frac{\partial_2 \phi(x - q)}{2})^2 + (\frac{D_2}{2} + \frac{\partial_1 \phi(x - q)}{2})^2 - \frac{\Delta \phi(x - q)}{4} - \lambda^2 = A^2 + B^2 - \frac{\Delta \phi(x-q)}{4} - \lambda^2 \\
    u_q(x) &:= u(x-q) \\
    \tilde{p}(x, \xi) &:= (\frac{\xi_1}{2} - \frac{\partial_2 \phi(x)}{2})^2 + (\frac{\xi_2}{2} + \frac{\partial_1 \phi(x)}{2})^2 + 1 = a(x, \xi)^2 + b(x, \xi)^2 + 1
\end{align*}

Note that, by rearranging \eqref{geigenfunction_eq}, we have

\[ \left((\frac{D_1}{2} - \frac{\partial_2 \phi}{2})^2 + (\frac{D_2}{2} + \frac{\partial_1 \phi}{2})^2 - \frac{\Delta \phi}{4} - \lambda^2 \right) u = Pu = 0\]

We begin with a few lemmas.

\begin{lemma} \thlabel{ABlemma}
    For $u$ such that $Pu = 0$, we have

    \begin{equation}
        ||Au||_{L^2(\R^2)} \lesssim ||u||_{L^2(\R^2)} \text{ and } ||Bu||_{L^2(\R^2)} \lesssim ||u||_{L^2(\R^2)}
    \end{equation}

    \begin{proof}

    The proof is similar to the proof of \thref{ABh_lemma} and is therefore omitted.
%
%
%
%
%
    \end{proof}
\end{lemma}

\begin{lemma} \thlabel{P_beta_lemma}
    For $\beta \in C_0^\infty(\R^2)$ and $u$ such that $Pu = 0$,

    \begin{equation}
        \sup_{q \in \R^2}||P(\beta_q u)||_{L^2(\R^2)} \lesssim ||u||_{L^2(\R^2)}
    \end{equation} 

\begin{proof}

The proof is similar to the proof of \thref{glcutoffefunc} and is therefore omitted. 

\end{proof}
    
\end{lemma}

\begin{lemma} \thlabel{Ae_beta_u_lemma}
    Suppose $\beta \in C_0^\infty(\R^2), \beta \equiv 1$ for $|x| < 1$, and $\beta \equiv 0$ for $|x| \geq 2$ and let $\beta_q(x) = \beta(x-q)$ for $q \in \R^2$. Then

    \begin{equation*}
        ||A(e^{i \sigma(\cdot, \nabla \phi(q))} \beta u_{-q})||_{L^2(\R^2)} \lesssim ||u||_{L^2(\R^2)} \text{ and } ||B(e^{i \sigma(\cdot, \nabla \phi(q))} \beta u_{-q})||_{L^2(\R^2)} \lesssim ||u||_{L^2(\R^2)}
    \end{equation*}

    \begin{proof}
        We only show the first inequality, as the proof for the second is nearly identical.

        \begin{align*}
            ||A(e^{i \sigma(\cdot, \nabla \phi(q))} \beta u_{-q})||_{L^2(\R^2)} &= ||e^{-i \sigma(\cdot, \nabla \phi(q))}A(e^{i \sigma(\cdot, \nabla \phi(q))} \beta u_{-q})||_{L^2(\R^2)} \\ &= ||(A_{-q} + \frac{\partial_2 \phi(\cdot+q)}{2} - \frac{\partial_2 \phi (q)}{2} - \frac{\partial_2 \phi(x)}{2}) \beta u_{-q}||_{L^2(\R^2)} \\ &\leq ||A_{-q} (\beta u_{-q})||_{L^2(\R^2)} + ||(\partial_2 \phi(\cdot + q) - \partial_2 \phi(q)) \beta u_{-q} ||_{L^2(\R^2)} + ||\partial_2 \phi \beta u_{-q}||_{L^2(\R^2)} \\ &\lesssim ||A (\beta_q u)||_{L^2(\R^2)} + ||\nabla \partial_2 \phi||_{L^\infty(\R^2)} ||u||_{L^2(\R^2)} + ||\partial_2 \phi||_{L^\infty(B(0,2))} ||u||_{L^2(\R^2)}
        \end{align*}

        Invoking \thref{ABlemma} and a simple estimate involving the commutator of $A$ and $\beta_q$ gives us the result. 
    \end{proof}
\end{lemma}

Now, we turn to the main theorem of this section, which is the following.

\begin{theorem}
    Suppose $\beta \in C_0^\infty(\R^2), \beta \equiv 1 \text{ for } |x| < 1, \text{ and } \beta \equiv 0 \text{ for } |x| \geq 2$ and let $\beta_q(x) = \beta(x - q)$ for $q \in \R^2$. Then

    \begin{equation}
        ||\beta_q u||_{L^\infty(\R^2)} \lesssim ||P (\beta_q u) ||_{L^2(\R^2)} + ||u||_{L^2(\R^2)}
    \end{equation}

    uniformly in $q$.

    \begin{proof}

    We have

    \begin{align*}
            ||\beta_q u||_{L^\infty(\R^2)} &= ||e^{i \sigma(\cdot, \nabla \phi(q))}\beta u_{-q}||_{L^\infty(\R^2)} \\ &\lesssim ||\mathcal{F}(e^{i \sigma(\cdot, \nabla \phi(q))} \beta u_{-q})||_{L^1(\R^2)} \\
            &\leq ||\frac{1}{\tilde{p}(0, \cdot)}||_{L^2(\R^2)} ||\tilde{p}(0, \cdot) \mathcal{F}(e^{i \sigma(\cdot, \nabla \phi(q))} \beta u_{-q})||_{L^2(\R^2)} \\ &\lesssim ||\mathcal{F}^{-1} \tilde{p}(0, \cdot) \mathcal{F}(e^{i \sigma(\cdot, \nabla \phi(q))} \beta u_{-q})||_{L^2(\R^2)}
    \end{align*}

    Note that $\mathcal{F}^{-1} \tilde{p}(0, \cdot) \mathcal{F}u = \tilde{p}^\mathrm{w}(0, D)u(x)$. Hence, we have

    \begin{align*}
        ||\beta_q u||_{L^\infty(\R^2)} &\lesssim ||\tilde{p}^\mathrm{w}(0, D) (e^{i \sigma(\cdot, \nabla \phi(q))} \beta u_{-q})||_{L^2(\R^2)} \\ &\leq ||(\tilde{p}^\mathrm{w}(0, D) - \tilde{p}^\mathrm{w}(\cdot, D))(e^{i \sigma(\cdot, \nabla \phi(q))} \beta u_{-q})||_{L^2(\R^2)} \\ &+ ||(\tilde{p}^\mathrm{w}(\cdot, D) - p^\mathrm{w}(\cdot, D))(e^{i \sigma(\cdot, \nabla \phi(q))} \beta u_{-q})||_{L^2(\R^2)} \\ &+ ||p^\mathrm{w}(\cdot, D)(e^{i \sigma(\cdot, \nabla \phi(q))} \beta u_{-q})||_{L^2(\R^2)}
    \end{align*}

    We will consider each of these summands in order. Denote

    \begin{align*}
    I &:= ||(\tilde{p}^\mathrm{w}(0, D) - \tilde{p}^\mathrm{w}(\cdot, D))(e^{i \sigma(\cdot, \nabla \phi(q))} \beta u_{-q})||_{L^2(\R^2)} \\
    II &:= ||(\tilde{p}^\mathrm{w}(\cdot, D) - p^\mathrm{w}(\cdot, D))(e^{i \sigma(\cdot, \nabla \phi(q))} \beta u_{-q})||_{L^2(\R^2)} \\
    III &:= ||p^\mathrm{w}(\cdot, D)(e^{i \sigma(\cdot, \nabla \phi(q))} \beta u_{-q})||_{L^2(\R^2)}
    \end{align*}
    
    We have

\begin{align} \label{p0minuspx}
    \tilde{p}^\mathrm{w}(0,D) - \tilde{p}^\mathrm{w}(x, D) = (\frac{D_1}{2} - \frac{\partial_2 \phi(0)}{2})^2 - (\frac{D_1}{2} - \frac{\partial_2 \phi(x)}{2})^2 + (\frac{D_2}{2} + \frac{\partial_1 \phi(0)}{2})^2 - (\frac{D_2}{2} + \frac{\partial_1 \phi(x)}{2})^2 
\end{align}

We will demonstrate that 

\begin{equation*}
    || ((\frac{D_1}{2} - \frac{\partial_2 \phi(0)}{2})^2 - (\frac{D_1}{2} - \frac{\partial_2 \phi(\cdot)}{2})^2) (e^{i \sigma(\cdot, \nabla \phi(q))} \beta u_{-q})||_{L^2(\R^2)} \lesssim ||u||_{L^2(\R^2)}
\end{equation*}

The steps to show a similar inequality for the remaining terms in \eqref{p0minuspx} are essentially the same. 

\begin{align*}
    &|| ((\frac{D_1}{2} - \frac{\partial_2 \phi(0)}{2})^2 - (\frac{D_1}{2} - \frac{\partial_2 \phi(\cdot)}{2})^2) (e^{i \sigma(\cdot, \nabla \phi(q))} \beta u_{-q})||_{L^2(\R^2)} \\ &\lesssim ||(\partial_2 \phi(0)^2 - \partial_2 \phi(\cdot)^2 )(e^{i \sigma(\cdot, \nabla \phi(q))} \beta u_{-q})||_{L^2(\R^2)} + ||(\partial_2 \phi(0) - \partial_2 \phi(\cdot))D_1(e^{i \sigma(\cdot, \nabla \phi(q))} \beta u_{-q})||_{L^2(\R^2)} \\ &+ ||(\partial_1 \partial_2 \phi)(e^{i \sigma(\cdot, \nabla \phi(q))} \beta u_{-q})||_{L^2(\R^2)}
\end{align*}

The first summand is bounded by $||(\partial_2 \phi)^2||_{L^\infty(B(0,2))}||u||_{L^2(\R^2)}$ because we are localized in $B(0,2)$. The third summand is bounded by $||\partial_1 \partial_2 \phi||_{L^\infty(B(0,2))}||u||_{L^2(\R^2)}$. Now, regarding the second summand, we have

\begin{align*}
    ||(\partial_2 \phi(0) - \partial_2 \phi(\cdot))D_1(e^{i \sigma(\cdot, \nabla \phi(q))} \beta u_{-q})||_{L^2(\R^2)} &\lesssim ||\partial_2 \phi||_{L^\infty(B(0,2))} ||D_1(e^{i \sigma(\cdot, \nabla \phi(q))} \beta u_{-q})||_{L^2(\R^2)} \\ &\lesssim ||A (e^{i \sigma(\cdot, \nabla \phi(q))} \beta u_{-q})||_{L^2(\R^2)} + ||\partial_2 \phi (e^{i \sigma(\cdot, \nabla \phi(q))} \beta u_{-q})||_{L^2(\R^2)} \\ &\lesssim ||u||_{L^2(\R^2)}
\end{align*}

where we have used \thref{Ae_beta_u_lemma}. Hence, $I \lesssim ||u||_{L^2(\R^2)}$.

    Because $\tilde{p}(x, \xi) - p(x, \xi) = 1 + \frac{\Delta \phi(x)}{4} + \lambda^2$, we have, by Theorem 4.3 in \cite{Zworski}, that

    \begin{equation*}
        (\tilde{p}^\mathrm{w}(x, D) - p^\mathrm{w}(x, D))u(x) = (1 + \frac{\Delta \phi(x)}{4} + \lambda^2)u
    \end{equation*}

    Thus,

    \begin{equation*}
        II = ||(\tilde{p}^\mathrm{w}(\cdot, D) - p^\mathrm{w}(\cdot, D))(e^{i \sigma(\cdot, \nabla \phi(q))} \beta u_{-q})||_{L^2(\R^2)} \leq (1 + \frac{||\Delta \phi||_{L^\infty(\R^2)}}{4} + \lambda_{\max}^2)||u||_{L^2(\R^2)}
    \end{equation*}

Finally, 

\begin{align*}
    III &= || e^{-i \sigma(\cdot, \nabla \phi(q))} P (e^{i \sigma(\cdot, \nabla \phi(q))} \beta u_{-q})||_{L^2(\R^2)} \\ &= ||(A_{-q} + \frac{\partial_2 \phi(\cdot+q)}{2} - \frac{\partial_2 \phi (q)}{2} - \frac{\partial_2 \phi(x)}{2})^2 + (B_{-q} - \frac{\partial_1 \phi(\cdot+q)}{2} + \frac{\partial_1 \phi (q)}{2} + \frac{\partial_1 \phi(x)}{2})^2 - \frac{\Delta \phi}{4} - \lambda^2) \beta u_{-q}||_{L^2(\R^2)} \\ &\lesssim ||P_{-q} \beta u_{-q}||_{L^2(\R^2)} + ||u||_{L^2(\R^2)}
\end{align*}

    \end{proof}
\end{theorem}

\section{Acknowledgments}

I'd like to thank my advisors, Matei Machedon and Manoussos Grillakis, for all the help they've provided me in bringing this work to fruition. I'd furthermore like to thank Xiaoqi Huang for early discussions we had that helped me at the beginning of this project.

\section{Declaration of Interest Statement}

I declare no conflict of interest.

\bibliography{bibliography}{}
\bibliographystyle{plain}

\section{Appendix}

We define the symbol classes

\begin{equation*}
    S_{\rho, \delta}(m) := \{ a \in C^\infty | |\partial_\xi^\alpha \partial_x^\beta a| \lesssim h^{-\rho |\alpha| - \delta |\beta|} m \text{ for all multiindices } \alpha, \beta \} 
\end{equation*}

which are an extension of Zworski's symbol classes in \cite{Zworski}. Note that $S_{\rho, \delta}(m) \subseteq S_{\max(\rho, \delta)}(m)$, so we  make use of many of the results that Zworski already proved. The important difference is that these symbol classes have better asymptotic expansions, and the important case for us is the following.

\begin{lemma}
    If $a \in S_{0, 1/2}(m_1)$ and $b \in S_{0, 1/2}(m_2)$, then

    \begin{equation*}
        a \# b = ab + \frac{h}{2i} \{a, b\} + O_{S_{1/2}(m_1m_2)}(h^{3/2})
    \end{equation*}
\end{lemma}

Importantly, the terms of the asymptotic expansion improve, i.e. have higher powers of $h$, which would not be the case if we simply worked with $S_{1/2}(m_1)$ and $S_{1/2}(m_2)$.

Below, we collect a few important results that we make use of above. See \cite{Zworski} for the proofs of the results below.

\begin{lemma} \thlabel{commutator_lemma}
    \begin{equation*}
        (D_{x_j} a)^\mathrm{w} = [D_{x_j}, a^\mathrm{w}] \text{ and } h(D_{\xi_j}a)^\mathrm{w} = - [x_j, a^\mathrm{w}]
    \end{equation*}

    for $j = 1, ..., n$.
\end{lemma}

\begin{lemma} ($L^2$ boundedness for symbols in $S$) \thlabel{L2_bound_for_S}

If the symbol $a$ belongs to $S_\delta$ for some $0 \leq \delta \leq 1/2$, then

\begin{equation*}
    a^\mathrm{w}(x, hD): L^2(\R^n) \to L^2(\R^n)
\end{equation*}

is bounded, with the estimate

\begin{equation*}
    ||a^\mathrm{w}(x, hD)||_{L^2 \to L^2} \leq C \sum_{|\alpha| \leq Mn} h^{|\alpha|/2} \sup_{\R^n} |\partial^\alpha a|
\end{equation*}
    
\end{lemma}

See \cite{ktz} for the proofs of the results below.

\begin{lemma} (Bernstein) \thlabel{Bernstein}
    Suppose there exists $\chi \in C_0^\infty(\R^k \times \R^k)$ and $N \geq 0$ such that

    \begin{equation*}
        u_h = \chi(x, hD) u_h + \mathcal{O}_{\mathcal{S}}(h^\infty), ||u_h||_{L^2} = \mathcal{O}(h^{-N})
    \end{equation*}

    Then, for any $1 \leq q \leq p \leq \infty$,

    \begin{equation*}
        ||u_h||_{L^p} \lesssim h^{k(1/p - 1/q)}||u_h||_{L^q} + \mathcal{O}(h^\infty)
    \end{equation*}
\end{lemma}

\begin{lemma} (Energy Estimate) \thlabel{energy_estimate}
    Suppose $a \in S(\R \times \R^k \times \R^k)$ is real-valued and that

    \begin{equation*}
        (hD_t + a^\mathrm{w}(t, x,hD_x))u(t,x)= f(t,x), u(0, x) = u_0(x)
    \end{equation*}
    \begin{equation*}
        f \in L^2(\R \times \R^k), u_0 \in L^2(\R^k)
    \end{equation*}

    Then

    \begin{equation*}
        ||u(t, \cdot)||_{L^2(\R^k)} \leq \frac{\sqrt{t}}{h} ||f||_{L^2(\R \times \R^k)} + ||u_0||_{L^2(\R^k)}
    \end{equation*}
\end{lemma}

\begin{lemma} (Strichartz) \thlabel{strichartz}
    Suppose $A =a^w(x, hD)$ and $\chi \in C_0^\infty(\R \times \R)$ is such that $\partial_\xi^2a$ is non-degenerate in the support of $\chi$. $F(t,r)$ is defined by $hD_t F(t,r) + A(t) F(t,r) = 0, F(r,r) = I$ and $\psi \in C_0^\infty(\R)$ with support sufficiently close to 0. Let

    \begin{equation*}
        U(t,r) = \psi(t) F(t,r) \chi^\mathrm{w} \text{ or } U(t,r) = \psi(t) \chi^\mathrm{w} F(t,r) 
    \end{equation*}

    Then
    \begin{equation*}
        ||\int_{t_1}^{t_2} U(t,s) f(s, x)ds ||_{L^6(\R^2)} \lesssim h^{-1/6} \int_{\R} ||f(s, x)||_{L^2(\R)}ds
    \end{equation*}
    
\end{lemma}

\end{document}